\documentclass[twoside,12pt]{amsart}
\usepackage{amsmath, amssymb, latexsym, amsthm, mathrsfs, color, geometry}
\usepackage[all]{xy}

\newcommand{\ed}{

\subsection*{Acknowledgments}
A part of this work was carried out during a visit of the first named author at the Kurt G\"odel Research Center, supported by FWF Grant M 1244-N13 of the second named author.
The first named author thanks the second named author for his kind hospitality,
and the Center's Director, researchers and staff for the excellent academic and friendly atmosphere.
We thank Arnold W. Miller for useful discussions.

\end{document}}

      \newenvironment{changemargin}[2]{\begin{list}{}{
         \setlength{\topsep}{0pt}\setlength{\leftmargin}{0pt}
         \setlength{\rightmargin}{0pt}
         \setlength{\listparindent}{\parindent}
         \setlength{\itemindent}{\parindent}
         \setlength{\parsep}{0pt plus 1pt}
         \addtolength{\leftmargin}{#1}\addtolength{\rightmargin}{#2}
         }\item }{\end{list}}

\newcommand{\nc}{\newcommand}
\nc{\productive}[2]{\bigl(#1,\allowbreak #2\bigr)^\x}
\nc{\set}[2]{{\{#1 \;:\; #2\}}}
\nc{\seq}[2]{{\la #1 \;:\; #2\ra}}
\nc{\cube}{(2^\w)^\w}
\nc{\Match}{\op{Match}}
\nc{\concat}[1]{\hat{\phantom{a}}\langle #1\rangle}
\nc{\poset}{\mathbb{P}}
\nc{\fn}[1]{{\op{Fn}(#1\times\w,2)}}
\nc{\linadd}{\op{linadd}}
\nc{\nonprod}{\non^\x}
\nc{\Ga}{\Gamma}
\nc{\Om}{\Omega}
\nc{\alephes}{{\aleph_0}}
\nc{\my}[1]{{\color{red} #1}}
\nc{\Cp}{\op{C}_p}
\nc{\Bp}{\op{B}_p}
\nc{\Pa}[8]{\bibitem{#1} {#2}, \emph{#3}, {#4} \textbf{#5} ({#6}), {#7}--{#8}.}
\nc{\tPa}[5]{\bibitem{#1} {#2}, \emph{#3}, {#4}, to appear.}
\nc{\sPa}[4]{\bibitem{#1} {#2}, \emph{#3}, {#4}, submitted.}
\nc{\Bc}[9]{\bibitem{#1} {#2}, \emph{#3}, in: \textbf{#4} (#5), #6 #7, #8--#9.}
\nc{\fD}{\mathfrak{D}}
\nc{\fX}{\mathfrak{X}}
\nc{\Onbd}{\Op_{\mathrm{nbd}}} 
\nc{\Omnb}{\Om_{\mathrm{nbd}}} 
\nc{\od}{\mathfrak{od}}
\nc{\Setting}[7]{\xymatrix@R=4pt@C=7pt{#1\ar@{-}[r]&#2\ar@{-}[r]&#3\\&#4\ar@{-}[u]\\
#5\ar@{-}[uu]\ar@{-}[r] & #6\ar@{-}[u]\ar@{-}[r] & #7\ar@{-}[uu]}}
\nc{\mx}[1]{\begin{matrix}#1\end{matrix}}
\nc{\plim}{p\txt{-}\lim}
\nc{\Bgp}{{\Z^\N}}
\nc{\Cgp}{{{\Z_2}^\N}}
\nc{\Cite}[1]{\textbf{[#1]}}
\nc{\Next}[1]{{#1^+}}
\nc{\Fr}{\mathit{F\!r}}
\nc{\intvl}[2]{{[#1(#2),\allowbreak #1(#2\!+\!1))}}
\nc{\Bdd}{\mathbf{B}}
\nc{\Ax}{\mathsf{Ax}}
\nc{\Dfin}{\mathfrak{D}_\mathrm{fin}}
\nc{\grbl}{{\mbox{\textit{\tiny gp}}}}
\nc{\bbP}{\mathbb{P}}
\nc{\BOfat}{\B_{\Om_{\mathrm{fat}}}}
\nc{\Bgood}{\B_{\mathrm{good}}}
\nc{\compactN}{\cl{\mathbb{N}}}
\nc{\blocks}[2]{\op{cl}_{#2}(#1)}
\nc{\blocksplus}[2]{\op{cl}^+_{#2}(#1)}
\nc{\arx}[1]{\texttt{http://arxiv.org/math/#1}}
\nc{\bq}{\begin{quote}}
\nc{\eq}{\end{quote}}
\nc{\cl}[1]{\overline{#1}}
\nc{\CH}{the Continuum Hypothesis}
\nc{\MA}{Martin's Axiom}
\nc{\Bfat}{\B_\mathrm{fat}}
\nc{\inv}{^{-1}}
\nc{\Cantor}{{2^\w}}
\nc{\bP}{\mathbf{P}}
\nc{\bof}{\op{\fb}}
\nc{\bofF}{\bof(\cF)}
\nc{\sr}[3]{{\txt{#1\\#2\\#3}}}
\nc{\gp}{\binom{\Om}{\Ga}}
\nc{\gpsmall}{\mbox{$\gp$}}
\nc{\gig}{\gimel}
\nc{\gns}{\sone(\Om,\gig)}
\nc{\nsr}[2]{#1}
\nc{\N}{\mathbb{N}}
\nc{\NN}{{\N^{\N}}}
\nc{\ZN}{{\Z^{\N}}}
\nc{\NNup}{{\N^{\uparrow\N}}}
\nc{\PN}{{P(\N)}}
\nc{\roth}{{[\w]^{\w}}}
\nc{\Fin}{{[\w]^{<\w}}}
\nc{\ici}{{[\w]^{(\w,\w)}}}
\nc{\Inc}{{\compactN^{\uparrow\N}}}
\nc{\powInc}[1]{{\big(\Inc\big)^{#1}}}
\nc{\powFin}[1]{{\big(\Fin\big)^{#1}}}
\nc{\powPN}[1]{{\big(\PN\big)^{#1}}}
\nc{\NcompactN}{{\compactN^\N}}
\nc{\setseq}[1]{\{#1 : n\in\N\}}
\nc{\sseq}[1]{\{#1 : n\in\N\}}
\nc{\Uarrow}{\smash{\big\uparrow}}
\nc{\LE}{\preccurlyeq}
\nc{\GE}{\succcurlyeq}
\nc{\op}{\operatorname}
\nc{\im}{\op{im}}
\nc{\Span}{\op{span}}
\nc{\maxfin}{\op{maxfin}}
\nc{\ran}{\op{range}}
\nc{\iso}{\cong}
\nc{\Madd}{{\M}^\star}
\nc{\cI}{\mathcal{I}}
\nc{\cJ}{\mathcal{J}}
\nc{\scrA}{\mathscr{A}}
\nc{\scrB}{\mathscr{B}}
\nc{\scrC}{\mathscr{C}}
\nc{\scrD}{\mathscr{D}}
\nc{\scrF}{\mathscr{F}}
\nc{\scrK}{\mathscr{K}}
\nc{\A}{\forall}
\nc{\B}{\mathrm{B}}
\nc{\cB}{\mathcal{B}}
\nc{\bB}{\mathbf{B}}
\nc{\BG}{\B_\Ga}
\nc{\BL}{\B_\Lambda}
\nc{\BT}{\B_\Tau}
\nc{\BTstar}{\B_{\Tau^*}}
\nc{\BO}{\B_\Om}
\nc{\DO}{\cD_\Om}
\nc{\KO}{\cK_\Om}
\nc{\CG}{C_\Ga}
\nc{\CL}{C_\Lambda}
\nc{\CT}{C_\Tau}
\nc{\CTstar}{C_{\Tau^*}}
\nc{\CO}{C_\Om}
\nc{\COgp}{C_{\Om^{\grbl}}}
\nc{\CLgp}{C_{\Lambda^{\grbl}}}
\nc{\BOgp}{\B_{\Om}^{\grbl}}
\nc{\BLgp}{\B_{\Lambda^{\grbl}}}
\nc{\sfC}{\mathsf{C}}
\nc{\sfD}{\mathsf{D}}
\nc{\bD}{\mathbf{D}}
\nc{\Tau}{\mathrm{T}}
\nc{\cA}{\mathcal{A}}
\nc{\cK}{\mathcal{K}}
\nc{\cD}{\mathcal{D}}
\nc{\cF}{\mathcal{F}}
\nc{\cS}{\mathcal{S}}
\nc{\cG}{\mathcal{G}}
\nc{\cY}{\mathcal{Y}}
\nc{\J}{\mathcal{J}}
\nc{\cL}{\mathcal{L}}
\nc{\cM}{\mathcal{M}}
\nc{\cN}{\mathcal{N}}
\nc{\cO}{\mathcal{O}}
\nc{\Op}{\mathrm{O}}
\nc{\cP}{\mathcal{P}}
\nc{\Q}{\mathbb{Q}}
\nc{\R}{\mathbb{R}}
\nc{\cU}{\mathcal{U}}
\nc{\Union}{\bigcup}
\nc{\cV}{\mathcal{V}}
\nc{\cW}{\mathcal{W}}
\nc{\Z}{{\mathbb Z}}
\nc{\Impl}{\Rightarrow}
\long\def\forget#1\forgotten{}
\nc{\ft}{\mathfrak{t}}
\nc{\fb}{\mathfrak{b}}
\nc{\fc}{\mathfrak{c}}
\nc{\fd}{\mathfrak{d}}
\nc{\fg}{\mathfrak{g}}
\nc{\oo}{\infty}
\nc{\fr}{\mathfrak{r}}
\nc{\fu}{\mathfrak{u}}
\nc{\fh}{\mathfrak{h}}
\nc{\fp}{\mathfrak{p}}
\nc{\fj}{\mathfrak{j}}
\nc{\fs}{\mathfrak{s}}
\nc{\w}{\omega}
\nc{\x}{\times}
\nc{\Iff}{\Leftrightarrow}

\nc{\nin}{\notin}
\nc{\cat}{\hat{\ }}
\nc{\sub}{\subseteq}
\nc{\spst}{\supseteq}
\nc{\sm}{\setminus}
\nc{\as}{\subseteq^*}
\nc{\rest}{\restriction}

\nc{\la}{\langle}
\nc{\ra}{\rangle}
\nc{\E}{\exists}
\nc{\dom}{\op{dom}}
\nc{\cov}{\op{cov}}
\nc{\add}{\op{add}}
\nc{\cof}{\op{cof}}
\nc{\cf}{\op{cf}}
\nc{\non}{\op{non}}
\nc{\unif}{\op{non}}
\nc{\COV}{\op{COV}}
\nc{\ADD}{\op{ADD}}
\nc{\COF}{\op{COF}}
\nc{\NON}{\op{NON}}
\nc{\impl}{\to}
\nc{\Lp}{\mathcal{L_\p}}
\nc{\Wlog}{without loss of generality}
\newtheorem{thm}{Theorem}[section]
\nc{\bthm}{\begin{thm}} \nc{\ethm}{\end{thm}}
\newtheorem{prop}[thm]{Proposition}
\nc{\bprp}{\begin{prop}} \nc{\eprp}{\end{prop}}
\newtheorem{fact}[thm]{Fact}
\nc{\bfct}{\begin{fact}} \nc{\efct}{\end{fact}}
\newtheorem{prob}[thm]{Problem}
\nc{\bprb}{\begin{prob}} \nc{\eprb}{\end{prob}}
\newtheorem{lem}[thm]{Lemma}
\nc{\blem}{\begin{lem}} \nc{\elem}{\end{lem}}
\newtheorem{claim}[thm]{Claim}
\nc{\bclm}{\begin{claim}} \nc{\eclm}{\end{claim}}
\newtheorem{cor}[thm]{Corollary}
\nc{\bcor}{\begin{cor}} \nc{\ecor}{\end{cor}}
\newtheorem{conj}[thm]{Conjecture}
\nc{\bcnj}{\begin{conj}} \nc{\ecnj}{\end{conj}}
\theoremstyle{definition}
\newtheorem{defn}[thm]{Definition}
\nc{\bdfn}{\begin{defn}} \nc{\edfn}{\end{defn}}
\theoremstyle{remark}
\newtheorem{rem}[thm]{Remark}
\nc{\brem}{\begin{rem}} \nc{\erem}{\end{rem}}
\newtheorem{cnv}[thm]{Convention}
\nc{\bcnv}{\begin{cnv}} \nc{\ecnv}{\end{cnv}}
\newtheorem{exam}[thm]{Example}
\nc{\bexm}{\begin{exam}} \nc{\eexm}{\end{exam}}
\nc{\bpf}{\begin{proof}} \nc{\epf}{\end{proof}}
\nc{\be}{\begin{enumerate}}
\nc{\ee}{\end{enumerate}}
\nc{\bi}{\begin{itemize}}
\nc{\itm}{\item}
\nc{\ei}{\end{itemize}}
\nc{\Subsection}[1]{\goodbreak\subsection*{#1}}

\forget
\forgotten

\nc{\sone}{\mathsf{S}_1}
\nc{\sfin}{\mathsf{S}_\mathrm{fin}}
\nc{\ufin}{\mathsf{U}_\mathrm{fin}}
\nc{\Split}{\mathsf{Split}}

\nc{\gone}{\mathsf{G}_1}    \nc{\gfin}{\mathsf{G}_\mathrm{fin}}

\title[The Gerlits--Nagy property and concentrated sets]{Additivity of the Gerlits--Nagy property and concentrated sets}
\author[B. Tsaban]{Boaz Tsaban}
\address[Tsaban]{Department of Mathematics, Bar-Ilan University, Ramat-Gan 52900, Israel}
\email{tsaban@math.biu.ac.il}
\urladdr{http://www.cs.biu.ac.il/\~{}tsaban}
\author[L. Zdomskyy]{Lyubomyr Zdomskyy}
\address[Zdomskyy]{Kurt G\"odel Research Center for Mathematical Logic, University of Vienna,
W\"ahringer Str.\ 25, 1090 Vienna, Austria}
\email{lzdomsky@logic.univie.ac.at}
\urladdr{http://www.logic.univie.ac.at/~lzdomsky/}

\subjclass{%
03E17 
26A03, 
03E75 
}

\begin{document}

\begin{abstract}
We settle all problems posed by Scheepers, in his tribute paper to Gerlits,
concerning the additivity of the Gerlits--Nagy property and related additivity numbers.
We apply these results to compute the minimal number of concentrated sets of reals
(in the sense of Besicovitch) whose union, when multiplied with a Gerlits--Nagy space,
need not have Rothberger's property. We apply these methods to construct a large family
of spaces, whose product with every Hurewicz space has Menger's property.
\end{abstract}

\maketitle

\section{Introduction}

We consider preservation of several classic topological properties under unions.
These properties are best understood
in the broader context of topological selection principles. We thus provide, in the present section,
a brief introduction.\footnote{This introduction
is adopted from \cite{LinSAdd}. Extended introductions
to this field are available in \cite{KocSurv, LecceSurvey, ict}.}
This framework was introduced by Scheepers in \cite{coc1}
to study, in a uniform manner, a variety of properties introduced
in different mathematical disciplines, since the early
1920's, by Menger, Hurewicz, Rothberger, Gerlits and Nagy, and many others.

By \emph{space} we mean an infinite topological space.
Let $X$ be a space. We say that $\cU$ is a \emph{cover}
of $X$ if $X=\Union\cU$, but $X\nin\cU$.
Often, $X$ is considered as a subspace of another space $Y$,
and in this case we always consider covers of $X$ by subsets of $Y$,
and require instead that no member of the cover contains $X$.
Let $\Op(X)$ be the family of all countable open covers of $X$.\footnote{Our assumption that
the considered covers are countable may be replaced by assuming that all considered spaces
are Lindel\"of in all finite powers, e.g., subsets of the real line.}
Define the following subfamilies of $\Op(X)$:
$\cU\in\Om(X)$ if each finite subset of $X$ is contained in some member of $\cU$.
$\cU\in\Ga(X)$ if $\cU$ is infinite, and each element of $X$ is contained in all but
finitely many members of $\cU$.

\begin{figure}[!htp]
\begin{changemargin}{-4cm}{-3cm}
\begin{center}
{
$\xymatrix@R=8pt{
&
&
& \sr{$\ufin(\Op,\Ga)$}{Hurewicz}{$\fb;\fb$}\ar[r]
& \sr{$\ufin(\Op,\Om)$}{$\fd; 2$}{}\ar[rr]
& & \sr{$\sfin(\Op,\Op)$}{Menger}{$\fd; \ge\max\{\fb,\fg\}$}
\\
&
&
& \sr{$\sfin(\Ga,\Om)$}{$\fd; 2$}{}\ar[ur]
\\
& \sr{$\sone(\Ga,\Ga)$}{$\fb; \ge\fh$}{}\ar[r]\ar[uurr]
& \sr{$\sone(\Ga,\Om)$}{$\fd; 2$}{}\ar[rr]\ar[ur]
& & \sr{$\sone(\Ga,\Op)$}{$\fd; \ge\add(\cN)$}{}\ar[uurr]
\\
&
&
& \sr{$\sfin(\Om,\Om)$}{Menger$^\uparrow$}{$\fd; 2$}\ar'[u][uu]
\\
& \sr{$\gp$}{Gerlits--Nagy}{$\fp; 2$}\ar[r]\ar[uu]
& \sr{$\sone(\Om,\Om)$}{Rothberger$^\uparrow$}{$\cov(\cM); 2$}\ar[uu]\ar[rr]\ar[ur]
& & \sr{$\sone(\Op,\Op)$}{Rothberger}{$\cov(\cM); \ge\add(\cN)$}\ar[uu]
}$
}
\caption{The Scheepers Diagram}\label{SchDiag}
\end{center}
\end{changemargin}
\end{figure}
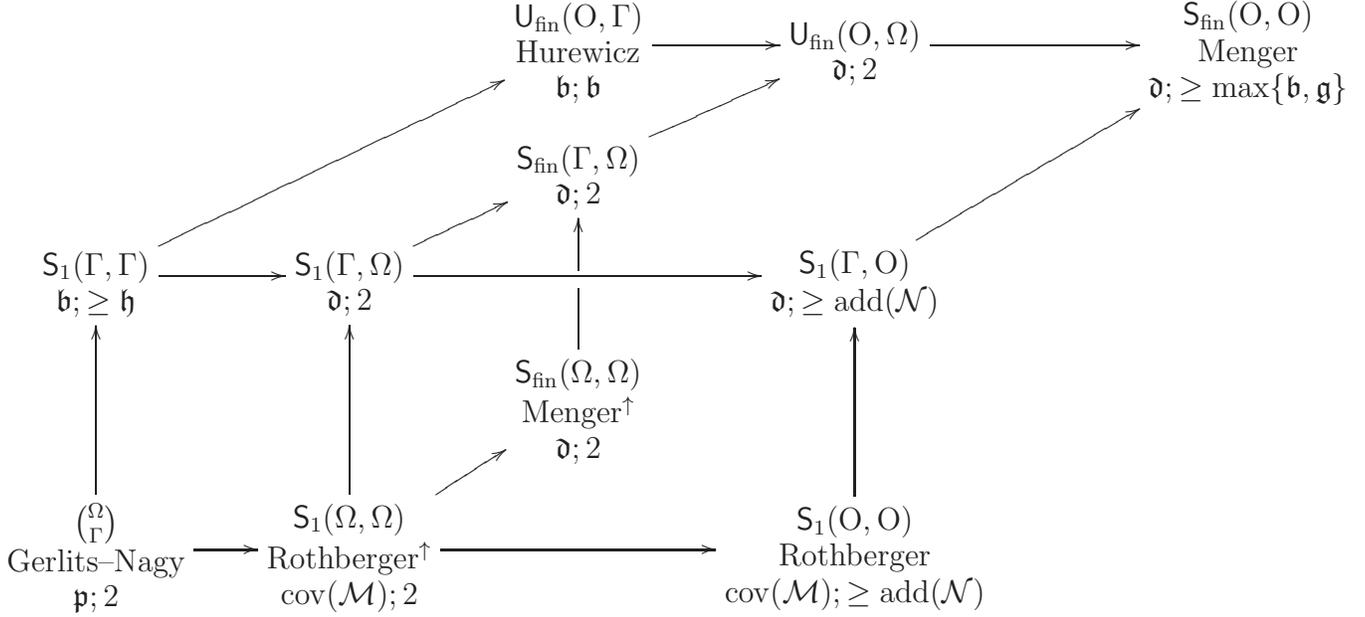

Some of the following statements may hold for families $\scrA$ and $\scrB$ of covers of $X$.
\begin{description}
\item[$\binom{\scrA}{\scrB}$] Each member of $\scrA$ contains a member of $\scrB$.
\item[$\sone(\scrA,\scrB)$] For each sequence $\seq{\cU_n\in\scrA}{n\in\N}$, there is a selection
$\seq{U_n\in\cU_n}{n\in\N}$ such that $\sseq{U_n}\in\scrB$.
\item[$\sfin(\scrA,\scrB)$] For each sequence $\seq{\cU_n\in\scrA}{n\in\N}$, there is a selection
of finite sets $\seq{\cF_n\sub\cU_n}{n\in\N}$ such that $\Union_n\cF_n\in\scrB$.
\item[$\ufin(\scrA,\scrB)$] For each sequence $\seq{\cU_n\in\scrA}{n\in\N}$, where no $\cU_n$ contains a finite
subcover, there is a selection
of finite sets $\seq{\cF_n\sub\cU_n}{n\in\N}$ such that $\sseq{\Union\cF_n}\in\scrB$.
\end{description}
We say, e.g., that $X$ satisfies $\sone(\Op,\Op)$ if the statement $\sone(\Op(X),\Op(X))$ holds.
This way, $\sone(\Op,\Op)$ is a property (or a class) of spaces, and similarly for all other statements
and families of covers.
Each nontrivial property among these properties, where $\scrA,\scrB$ range over $\Op,\Om,\Ga$,
is equivalent to one in Figure \ref{SchDiag}  \cite{coc1, coc2}. In this diagram, an arrow denotes implication.

The names indicated below some of the properties are of those introducing it.
The two names ending with a symbol $\uparrow$ indicate that the properties $\sone(\Om,\Om)$ and
$\sfin(\Om,\Om)$ are characterized by having all finite
powers Rothberger and Menger, respectively \cite{Sakai88, coc2}.
In addition, we indicate below each class $P$ its \emph{critical cardinality} $\non(P)$
(the minimal cardinality of a space not in the class), followed by its \emph{additivity number} $\add(P)$
(the minimal number of spaces in the class with union outside the class). When only upper and lower
bounds are known, we write a lower bound. To save space, we do not write the immediate upper bound,
$\cf(\non(P))$.
These cardinals are all combinatorial cardinal characteristics of the continuum,
details about which are available in \cite{BlassHBK}. Here,
$\cM,\cN$ are the families of meager sets in $\R$ and Lebesgue null sets in $\R$,
respectively.
Complete computations of the mentioned additivity numbers and bounds,
with references, are available in \cite{AddQuad}.
That the additivity
number of $\sone(\Ga,\Op)$ is $\ge\add(\cN)$ follows from Bartoszy\'nski's
Theorem \cite[Lemma 2.16]{AddQuad} and the first observation in \cite[Appendix A]{MHP}.

Many additional---classic and new---properties are studied in relation to the the Scheepers Diagram.
One of these is the Gerlits--Nagy property, to which we now focus our attention.

\section{Additivity of the Gerlits--Nagy property}

\bdfn
For classes $P,Q$ of spaces, $\add(P,Q)$ is the minimal cardinal $\kappa$ such that some
union of $\kappa$ members of $P$ is not in $Q$. $\add(P)$ is $\add(P,P)$.
\edfn

A countable cover $\cU$ of a space $X$ is in $\gimel(\Gamma)$ ($\gig$, read \emph{gimel}, for brevity)\footnote{In
general, the gimel operator $\gimel$ can be applied to any type of covers. However, in the present
paper we apply it only to $\Gamma$.}
if for each (equivalently, some)
bijective enumeration $\cU=\sseq{U_n}$, there is an increasing $h\in\NN$ such that, for each $x\in X$,
$$x\in \Union_{k=h(n)}^{h(n+1)-1}U_k$$
for all but finitely many $n$.

The property $\gns$ was introduced, in an equivalent form, by Gerlits and Nagy in their seminal paper \cite{GN}.
Building on results of Gerlits and Nagy and extending them, Ko\v{c}inac and Scheepers prove in \cite{coc7}
that
$$\ufin(\Op,\Ga)\cap\sone(\Op,\Op)=\gns.$$
This property is often referred to as the \emph{Gerlits--Nagy property} \cite{SchGerlits}.

The importance of the Gerlits--Nagy property $\gns$
in various contexts is surveyed in Scheepers's tribute to Gerlits \cite{SchGerlits}.
In \cite[\S{} II.5]{SchGerlits}, Scheepers poses several problems concerning preservation of this
property under unions. Scheepers's problems are all settled by the following two theorems.

\bthm\label{yes}
$$\add(\gns,\sone(\Op,\Op))=\add(\gpsmall,\sone(\Op,\Op))=\cov(\cM).$$
\ethm
\bpf
Since $\gp=\sone(\Om,\Ga)$ \cite{GN} implies $\gns$,
$$
\add(\gns,\sone(\Op,\Op))\le\add(\gpsmall,\sone(\Op,\Op))\le
\non(\sone(\Op,\Op))=\cov(\cM).
$$
It remains to prove that $\cov(\cM)\le\add(\gns,\sone(\Op,\Op))$.
We use that $\sone(\Op,\Op)=\sone(\Om,\Op)$ \cite{coc1}.

Let $\kappa<\cov(\cM)$.
Assume that, for each $\alpha<\kappa$, $X_\alpha$ satisfies $\gns$, and
$X=\Union_{\alpha<\kappa}X_\alpha$.
Let $\cU_n\in\Om(X)$ for all $n$. Enumerate $\cU_n=\set{U^n_m}{m\in\N}$.
For each $\alpha$, as $X_\alpha$ satisfies $\gns$,
there are $f_\alpha\in\NN$ and an increasing $h_\alpha\in\NN$ such that,
for each $x\in X_\alpha$,
$$x\in \Union_{k=h_\alpha(n)}^{h_\alpha(n+1)-1}U^k_{f_\alpha(k)}$$
for all but finitely many $n$.

Since $\kappa<\cov(\cM)\le\fd$ \cite{BlassHBK}, there is an increasing $h\in\NN$
such that, for each $\alpha<\kappa$, the set
$$I_\alpha = \set{n}{[h_\alpha(n),h_\alpha(n+1))\sub [h(n),h(n+1))}$$
is infinite \cite{BlassHBK}.
For each $\alpha<\kappa$, define
$$g_\alpha\in \prod_{n\in I_\alpha} \N^{[h(n),h(n+1))}$$
by
$g_\alpha(n)=f_\alpha\rest[h(n),h(n+1))$ for all $n\in I_\alpha$.
As $\kappa<\cov(\cM)$, by Lemma 2.4.2(3) in \cite{BarJu},
there is $g\in\prod_n \N^{[h(n),h(n+1))}$ guessing all functions $g_\alpha$,
that is, for each $\alpha<\kappa$, $g(n)=g_\alpha(n)$ for infinitely many $n\in I_\alpha$ \cite{BlassHBK}.
Define $f\in\NN$ by $f(k)=g(n)(k)$, where $n$ is the one with $k\in[h(n),h(n+1))$.
Then $U^n_{f(n)}\in\Op(X)$.

Indeed, let $x\in X$. Pick $\alpha<\kappa$ with $x\in X_\alpha$.
Pick $m$ such that, for all $n>m$,
$x\in \Union_{k=h_\alpha(n)}^{h_\alpha(n+1)-1}U^k_{f_\alpha(k)}$.
Pick $n\in I_\alpha$ such that $n>m$ and $g(n)=g_\alpha(n)$. Then
$$x\in \Union_{k=h_\alpha(n)}^{h_\alpha(n+1)-1}U^k_{f_\alpha(k)}\sub
\Union_{k=h(n)}^{h(n+1)-1}U^k_{f_\alpha(k)}=\Union_{k=h(n)}^{h(n+1)-1}U^k_{f(k)}.\qedhere$$
\epf

We can now compute the additivity number of the Gerlits--Nagy property.

\bthm\label{yes2}
$\add(\gns)=\add(\cM)$.
\ethm
\bpf
As $\gns=\ufin(\Op,\Ga)\cap \sone(\Op,\Op)$,
\begin{eqnarray*}
\add(\gns) & \le & \non(\gns)=\\
& = & \min\{\non(\ufin(\Op,\Ga)),\non(\sone(\Op,\Op))\} =\\
& = & \min\{\fb,\cov(\cM)\}=\add(\cM).
\end{eqnarray*}
It remains to prove the other inequality.
Let $X=\Union_{\alpha<\kappa}X_\alpha$, with each $X_\alpha$ in $\gns$,
and $\alpha<\add(\cM)$. By Theorem \ref{yes}, $X$ satisfies $\sone(\Op,\Op)$.
As $\kappa<\fb=\add(\ufin(\Op,\Ga))$ \cite{AddQuad}, $X$ satisfies $\ufin(\Op,\Ga)$, too.
Thus, $X$ satisfies $\gns$.
\epf

The following definition and corollary will be used in the next section.

\bdfn
Let $P,Q$ be classes of spaces, each containing all one-element spaces and closed under
homeomorphic images.
$\productive{P}{Q}$ is the class of all spaces $X$ such that, for each $Y$ in $P$,
$X\x Y$ is in $Q$. $\productive{P}{P}$ is denoted $P^\x$.
\edfn

\blem\label{bds}
Let $P,Q$ be classes of spaces. Then:
\be
\itm $\add(P,Q)\le\non(\productive{P}{Q})\le\non(Q)$.
\itm $\add(Q)\le\add(\productive{P}{Q})\le\non(Q)$. \qed
\ee
\elem

\bcor\label{help}
$$\non(\productive{\gns}{\sone(\Op,\Op)})=\non(\productive{\gpsmall}{\sone(\Op,\Op)})=\cov(\cM).$$
\ecor
\bpf
Theorem \ref{yes} and Lemma \ref{bds}.
\epf

\brem
The above proofs, verbatim, show that the results of the present section also apply in the case
where countable \emph{Borel} covers are considered instead of countable open covers.
\erem

\section{Unions of concentrated sets}

According to Besicovitch \cite{Besicovitch1, Besicovitch2},
a space $X$ is \emph{concentrated} if there is a countable
$D\sub X$ such that for each open $U\spst D$, $X\sm U$ is countable.
More generally, for a cardinal $\kappa$, a space $X$ is \emph{$\kappa$-concentrated}
if there is a countable $D\sub X$ such that for each open $U\spst D$, $|X\sm U|<\kappa$.
The classic examples of concentrated spaces are Luzin sets.
Modern examples are constructed from scales (e.g., \cite{sfh}).

Babinkostova and Scheepers proved in \cite{BabSch} that every concentrated metric space
belongs to $\productive{\gns}{\sone(\Op,\Op)}$.
In other words, for each concentrated metric space $C$,
if $Y$ satisfies $\ufin(\Op,\Ga)$ and $\sone(\Op,\Op)$, then $C\x Y$ satisfies $\sone(\Op,\Op)$.
We generalize this result in several ways.

\bthm\label{conc1}
\mbox{}
\be
\itm Let $\lambda$ be a regular uncountable cardinal $\le\cov(\cM)$.
The minimal number of $\lambda$-concentrated spaces, whose union is a regular space not satisfying
$\productive{\gns}{\sone(\Op,\Op)}$, is $\cov(\cM)$.
\itm The minimal number of $\cov(\cM)$-concentrated spaces, whose union is a regular space not satisfying
$\productive{\gns}{\sone(\Op,\Op)}$, is $\cf(\cov(\cM))$.
\ee
\ethm
\bpf
We prove both statements simultaneously.

There is a set of real numbers, of cardinality $\cov(\cM)$, that does not satisfy $\sone(\Op,\Op)$ \cite{coc2}.
Thus, the minimal number sought for is at most $\cov(\cM)$ for (1) and at most
$\cf(\cov(\cM))$ for (2).

Let $\lambda$ be a regular cardinal $\le\cov(\cM)$ for (1), and $\cov(\cM)$ for (2).
Let $\kappa<\cov(\cM)$ for (1), and $<\cf(\cov(\cM))$ for (2).

Let $C=\Union_{\alpha<\kappa}C_\alpha$ be a regular space, with each
$C_\alpha$ $\lambda$-concentrated on some countable set $D_\alpha\sub C_\alpha$.
Let $Y$ be a space satisfying $\gns$. We must prove that $C\x Y$ satisfies $\sone(\Op,\Op)$.

Let $K$ be a compact space containing $C$ as a subspace.
Let $\cU_n$, $n\in\N$, be countable covers of $C\x Y$ by sets open in $K\x Y$.
Let $D=\Union_{\alpha<\kappa}D_\alpha$. As $|D|=\kappa<\cov(\cM)$, we have by Corollary \ref{help}
that $D\x Y$ satisfies $\sone(\Op,\Op)$. Thus, pick $U_n\in\cU_n$, $n\in\N$, such that
$D\x Y\sub U:=\Union_nU_n$.

The Hurewicz property $\ufin(\Op,\Ga)$ is preserved by products with compact spaces, moving to
closed subspaces,
and continuous images \cite{coc2}.
Since $Y$ satisfies $\ufin(\Op,\Ga)$ and $K$ is compact, $K\x Y$ satisfies $\ufin(\Op,\Ga)$.
Thus, so does $K\x Y\sm U$. It follows that the projection $H$ of $K\x Y\sm U$
on the first coordinate, satisfies $\ufin(\Op,\Ga)$.
Note that
$$(K\sm H)\x Y\sub U.$$
The argument in the proof of \cite[Theorem 5.7]{coc2} generalizes to
regular spaces, to show that for $H,F$ disjoint subspaces of a regular space $K$
with $H$ $\ufin(\Op,\Ga)$, and $F$ F$_\sigma$, there is a G$_\delta$ set $G\sub K$
such that $G\spst F$ and $H\cap G=\emptyset$.

For each $\alpha<\kappa$, let $G_\alpha$ be a G$_\delta$ subset of $K$ such that
$D_\alpha\sub G_\alpha$ and $H\cap G_\alpha=\emptyset$. As $C_\alpha$ is $\lambda$-concentrated on
$D_\alpha$, $C_\alpha\sm G_\alpha$ is a countable union of sets of cardinality $<\lambda$.

As $\lambda$ has uncountable cofinality, $|C_\alpha\sm G_\alpha|<\lambda$.
Then
$$C\cap H\sub C\sm \bigcup_{\alpha<\kappa} G_\alpha\sub \Union_{\alpha<\kappa}C_\alpha\sm G_\alpha.$$
By splitting to cases $\lambda<\cov(\cM)$ and $\lambda=\cov(\cM)$, one sees that
$|C\cap H|<\cov(\cM)$ in both scenarios (1) and (2).
Thus, by Corollary \ref{help} again, $(C\cap H)\x Y$ satisfies $\sone(\Op,\Op)$, and there are $V_n\in\cU_n$,
$n\in\N$, such that $(C\cap H)\x Y\sub \Union_n V_n$. In summary,
$$C\x Y\sub ((K\sm H)\x Y) \cup ((C\cap H)\x Y)\sub \Union_{n\in\N}(U_n\cup V_n).$$
We have picked two sets (instead of one) from each cover $\cU_n$, but this is fine \cite{GFTM95}
(cf.\ \cite[Appendix A]{MHP}).
\epf

\bdfn
Let $\kappa$ be an infinite cardinal number.
Let $\scrC_0(\kappa)$ be the family of regular spaces of cardinality $<\kappa$. For successor ordinals $\alpha+1$,
let $C\in\scrC_{\alpha+1}(\kappa)$ if $C$ is regular, and:
\be
\itm either there is a countable $D\sub C$ with $C\sm U\in\scrC_\alpha(\kappa)$ for all open $U\spst D$;
\itm or $C$ is a union of less than $\cf(\kappa)$ members of $\scrC_\alpha(\kappa)$.
\ee
For limit ordinals $\alpha$, let $\scrC_\alpha(\kappa)=\Union_{\beta<\alpha}\scrC_\beta(\kappa)$.
\edfn

Babinkostova and Scheepers prove, essentially, that every member of $\scrC_{\alephes}(2)$ is in
$\productive{\gns}{\sone(\Op,\Op)}$ \cite{BabSch}. We use our methods to prove the following, stronger result.

For the following theorem, we recall from the Scheepers Diagram that $\add(\cN)\le\add(\sone(\Op,\Op))\allowbreak\le\cf(\cov(\cM))$.

\bthm\label{last}
The product of each member of $\scrC_{\add(\cN)}(\cov(\cM))$ with every member of $\gns$ satisfies $\sone(\Op,\Op)$.
\ethm
\bpf
We prove the stronger assertion, with $\add(\sone(\Op,\Op))$ instead of $\add(\cN)$.

For brevity, let $\scrC_{\alpha}:=\scrC_{\alpha}(\cov(\cM))$ for all $\alpha$.
By induction on $\alpha\le\add(\sone(\Op,\Op))$, we prove that
$$\scrC_{\alpha}\sub\productive{\gns}{\sone(\Op,\Op)}.$$
The proof is similar to that of Theorem \ref{conc1}, so we omit some of the explanations.

The case $\alpha=0$ is treated in Corollary \ref{help}.
For limit $\alpha$, there is nothing to prove.

$\alpha+1$: Let $C\in \scrC_{\alpha+1}$.
Let $K$ be a compact space containing $C$ as a subspace.
Let $Y$ be a space satisfying $\gns$.

First case: There is a countable $D\sub C$ with $C\sm U\in\scrC_\alpha$ for all open $U\spst D$.
Given $\cU_n\in\Op(C\x Y)$, pick $U_n\in\cU_n$, $n\in\N$, such that $D\x Y\sub U:=\Union_nU_n$.
Let $H$ be the projection of $K\x Y\sm U$ on the first coordinate.
Let $G$ be a G$_\delta$ subset of $K$ such that
$D\sub G$ and $H\cap G=\emptyset$. $C\sm G$ is a countable union of elements of $\scrC_\alpha$.
By the induction hypothesis and Corollary \ref{help}, $C\sm G\in \productive{\gns}{\sone(\Op,\Op)}$.
Then
$(C\sm G)\x Y$ satisfies $\sone(\Op,\Op)$, and there are $V_n\in\cU_n$,
$n\in\N$, such that $(C\cap H)\x Y\sub (C\sm G)\x Y\sub \Union_n V_n$. In summary,
$$C\x Y\sub ((K\sm H)\x Y) \cup ((C\cap H)\x Y)\sub \Union_{n\in\N}(U_n\cup V_n).$$

Second case: There are $\kappa<\cf(\cov(\cM))$ and $C_\beta\in \scrC_\alpha$, $\beta<\kappa$,
such that $C=\Union_{\beta<\kappa}C_\beta$. For each $\beta<\kappa$ with $C_\beta$ a union of less than
$\cf(\cov(\cM))$ members of
$$\scrC_{<\alpha}:=\Union_{\gamma<\alpha}\scrC_{\gamma},$$
we may take all elements in all of these unions instead of the original $C_\beta$'s.
Thus, we may assume that for each $C_\beta$ there is a countable (possibly empty) $D_\beta\sub C_\beta$ with
$$C_\beta\sm U\in \scrC_{<\alpha}$$
for all open $U\spst D_\beta$.
Let $D=\Union_{\beta<\kappa}D_\beta$. Then $|D|<\cov(\cM)$.

Given $\cU_n\in\Op(C\x Y)$, pick $U_n\in\cU_n$, $n\in\N$, such that $D\x Y\sub U:=\Union_nU_n$.
Let $H$ be the projection of $K\x Y\sm U$ on the first coordinate.
For each $\beta<\kappa$, let $G_\beta$ be a G$_\delta$ subset of $K$ such that
$D_\beta\sub G_\beta$ and $H\cap G_\beta=\emptyset$. Let $G=\Union_{\beta<\kappa}G_\beta$.
Now,
$$C\cap H\sub \Union_{\beta<\kappa}C_\beta\sm G_\beta,$$
where each $C_\beta\sm G_\beta$ is a countable union of elements of
$\scrC_{<\alpha}$. All in all, we arrive at a union of a family $\cF\sub \scrC_{<\alpha}$
with $|\cF|<\cf(\cov(\cM))$, and we must show that $\Union\cF\in \productive{\gns}{\sone(\Op,\Op)}$.
Indeed, for each $\gamma<\alpha$,
$$X_\gamma:=\Union(\cF\cap\scrC_{\gamma})\in \scrC_{\gamma+1}\sub \scrC_{\alpha}.$$
By the induction hypothesis, $X_\gamma\in \productive{\gns}{\sone(\Op,\Op)}$. By Corollary \ref{help}, since $\alpha<\add(\sone(\Op,\Op))$,
$$\Union\cF=\Union_{\gamma<\alpha}X_\gamma\in \productive{\gns}{\sone(\Op,\Op)}.$$
It follows that $(C\cap H)\x Y\sub (C\sm G)\x Y\sub \Union_n V_n$ for some $V_n\in\cU_n$, $n\in\N$,
and $C\x Y\sub \Union_{n\in\N}(U_n\cup V_n)$.
\epf

\section{Spaces whose product with Hurewicz spaces are Menger}

A space is \emph{$\sigma$-compact} if it is a union of countably many compact spaces.

\bdfn
For a cardinal $\lambda$, $K_\lambda$ is the family of all spaces that are unions of
less than $\lambda$ compact spaces. A space $X$ is \emph{$K_\lambda$-concentreted} if
there is a $\sigma$-compact subset $D\sub X$ such that $X\sm U\in K_\lambda$ for each open
$U\spst D$.
\edfn

Babinkostova and Scheepers proved in \cite{BabSch} that, for each concentrated metric space $C$,
if $Y$ has Hurewicz's property $\ufin(\Op,\Ga)$, then $C\x Y$ has Menger's property $\sfin(\Op,\Op)$.
We use the methods of the previous section to generalize this result.
Since the proofs are almost literal repetition of the corresponding ones in the previous section,
we omit some of the details.

The following is immediate from the definitions.

\blem[Folklore]\label{folk}
$\add(\ufin(\Op,\Ga),\sfin(\Op,\Op))=\fd$.
\elem

\blem\label{yey}
$K_\fd\sub \productive{\ufin(\Op}{\Ga),\sfin(\Op,\Op)}$.
\elem
\bpf
Each compact space is in $\ufin(\Op,\Ga)^\x$. Apply Lemma \ref{folk}.
\epf

\bthm\label{conc2}
\mbox{}
\be
\itm Let $\lambda$ be a regular uncountable cardinal $\le\fd$.
The minimal number of $K_\lambda$-concentrated spaces, whose union is a regular space not satisfying
$\productive{\ufin(\Op}{\Ga),\sfin(\Op,\Op)}$, is $\fd$.
\itm The minimal number of $K_\fd$-concentrated spaces, whose union is a regular space not satisfying
$\productive{\ufin(\Op}{\Ga),\sfin(\Op,\Op)}$, is $\cf(\fd)$.
\ee
\ethm
\bpf
There is a set of real numbers, of cardinality $\fd$, that does not satisfy $\sfin(\Op,\Op)$ \cite{coc2}.
Thus, the minimal number sought for is at most $\fd$ for (1) and at most $\cf(\fd)$ for (2).

Let $\lambda$ be a regular cardinal $\le\fd$ for (1), and $\fd$ for (2).
Let $\kappa<\fd$ for (1), and $<\cf(\fd)$ for (2).

Let $C=\Union_{\alpha<\kappa}C_\alpha$ be a regular space, with each
$C_\alpha$ $K_\lambda$-concentrated on some $\sigma$-compact set $D_\alpha\sub C_\alpha$.
Let $Y$ be a space satisfying $\ufin(\Op,\Ga)$. We must prove that $C\x Y$ satisfies $\sfin(\Op,\Op)$.

Let $K$ be a compact space containing $C$ as a subspace.
Let $\cU_n$, $n\in\N$, be countable covers of $C\x Y$ by sets open in $K\x Y$.
Let $D=\Union_{\alpha<\kappa}D_\alpha$. As $D\in K_\fd$, we have by Lemma \ref{yey}
that $D\x Y$ satisfies $\sfin(\Op,\Op)$. Thus, pick finite $\cF_n\in\cU_n$, $n\in\N$, such that
$D\x Y\sub U:=\Union_n\Union\cF_n$.

Since $Y$ satisfies $\ufin(\Op,\Ga)$ and $K$ is compact, the projection $H$ of $K\x Y\sm U$
on the first coordinate satisfies $\ufin(\Op,\Ga)$.
Note that
$$(K\sm H)\x Y\sub U.$$
For each $\alpha<\kappa$, let $G_\alpha$ be a G$_\delta$ subset of $K$ such that
$D_\alpha\sub G_\alpha$ and $H\cap G_\alpha=\emptyset$. As $C_\alpha$ is $K_\lambda$-concentrated on
$D_\alpha$, $C_\alpha\sm G_\alpha$ is a countable union of elements of $K_\lambda$.

As $\lambda$ has uncountable cofinality, $C_\alpha\sm G_\alpha\in K_\lambda$.
Then
$$C\cap H\sub C\sm \bigcup_{\alpha<\kappa} G_\alpha\sub \tilde C:=\Union_{\alpha<\kappa}C_\alpha\sm G_\alpha.$$
By splitting to cases $\lambda<\fd$ and $\lambda=\fd$, one sees that
$\tilde C\in K_\fd$ in both scenarios (1) and (2).
Thus, by Lemma \ref{yey} again, $\tilde C\x Y$ satisfies $\sfin(\Op,\Op)$, and there are finite $\tilde \cF_n\in\cU_n$,
$n\in\N$, such that $\tilde C\x Y\sub \Union_n \Union\tilde\cF_n$.
Thus,
$$C\x Y\sub ((K\sm H)\x Y) \cup ((C\cap H)\x Y)\sub ((K\sm H)\x Y) \cup (\tilde C\x Y)\sub
\Union_{n\in\N}\Union(\cF_n\cup \tilde\cF_n).\qedhere$$
\epf

\bdfn
Let $\kappa$ be an infinite cardinal number.
Let $\scrK_0(\kappa)$ be the family of regular spaces in $K_\kappa$.
For successor ordinals $\alpha+1$,
let $C\in\scrK_{\alpha+1}(\kappa)$ if $C$ is regular, and:
\be
\itm either there is a $\sigma$-compact $D\sub C$ with $C\sm U\in\scrK_\alpha(\kappa)$ for all open $U\spst D$;
\itm or $C$ is a union of less than $\cf(\kappa)$ members of $\scrK_\alpha(\kappa)$.
\ee
For limit ordinals $\alpha$, let $\scrK_\alpha(\kappa)=\Union_{\beta<\alpha}\scrK_\beta$.
\edfn

For every $\alpha$ the class $\mathcal K_\alpha(\kappa)$ is closed under
products with compact regular spaces. In particular,
the classes $\scrK_\alpha(\kappa)$ are much wider than $\scrC_\alpha(\kappa)$.
Babinkostova and Scheepers prove, essentially, that every member of $\scrC_{\alephes}(2)$ is in
$\productive{\ufin(\Op}{\Ga),\sfin(\Op,\Op)}$ \cite{BabSch}.

\bthm
The product of each member of $\scrK_{\max\{\fb,\fg\}}(\fd)$ with every member of $\ufin(\Op,\Ga)$ satisfies $\sfin(\Op,\Op)$.
\ethm
\bpf
We recall from the Scheepers Diagram that $\max\{\fb,\fg\}\le\add(\sfin(\Op,\Op))\le\cf(\fd)$.
A combination of the arguments in the proofs of Theorems \ref{conc2} and \ref{last}
show that
$$\scrK_{\add(\sfin(\Op,\Op))}(\fd)\sub \productive{\ufin(\Op}{\Ga),\sfin(\Op,\Op)}.$$
Since we have already presented three proofs using these methods, we leave the verification to the
reader.
\epf

\ed